\documentclass[a4paper,12pt]{amsart} \usepackage{amssymb,euscript,eufrak,verbatim}
\usepackage[latin1]{inputenc}      
\usepackage[final]{epsfig}
\usepackage{psfrag}

\usepackage{tikz}
\usetikzlibrary{snakes}
\usetikzlibrary{arrows}

\def\e{\varepsilon} 
\def\L{\mathcal{L}}
\def\Z{\mathbb{Z}} 
\def\R{\mathbb{R}} 
\def\Si{\Sigma_2} 
\def\N{\mathbb{N}}

\def\Ts{T_{\hat{s}}}

\def\L{\mathcal{L}}
\def\D{\mathcal D}
\def\S1{\mathbb{S}^1}

\newtheorem{theorem}{Theorem}

\DeclareMathSymbol{\varnothing}{\mathord}{AMSb}{"3F} 
\begin{document}
\renewenvironment{proof}{\noindent {\bf Proof.}}{ \hfill\qed\\ }
\newenvironment{proofof}[1]{\noindent {\bf Proof of #1.}}{ \hfill\qed\\ }

\title{Recurrence in generic staircases}  
\author{Serge Troubetzkoy}
\address{Centre de physique th\'eorique\\
Federation de Recherches des Unites de Mathematique de Marseille\\
Institut de math\'ematiques de Luminy and\\ 
Universit\'e de la M\'editerran\'ee\\ 
Luminy, Case 907, F-13288 Marseille Cedex 9, France}
\email{troubetz@iml.univ-mrs.fr}
\urladdr{http://iml.univ-mrs.fr/{\lower.7ex\hbox{\~{}}}troubetz/} \date{}
\begin{abstract}  
The straight-line flow on almost every staircase and on almost every square tiled staircase is recurrent.
For almost every square tiled staircase the set of periodic orbits is dense in the phase space.
\end{abstract} 
\maketitle
\section{Introduction}
A compact translation surface is a surface which can be obtained by edge-to-edge gluing of finitely many polygons in the plane
using only translations.  Since the seminal work of Veech in 1989 \cite{Ve} the study of compact translation surfaces of finite area have developed extensively. The
study of translation surfaces of infinite area, obtained by gluing countably many polygons via translations, has only recently begun. 
A natural class of infinite translation
surfaces, staircases, were introduced in \cite{HuWe} and studied in \cite{HoWe}.
Billiards in irrational polygons give rise to another class of infinite translation surfaces \cite{GuTr}. 

One of the first dynamic properties of infinite translation surfaces one needs to understand
is the almost sure recurrence of the straight-line flow.  
Recurrence of infinite translation surfaces have been investigated in
\cite{GuTr}, \cite{Ho}, \cite{HoWe}, \cite{HuWe},  \cite{HuLeTr}, \cite{ScTr}, and \cite{Tr}. Hubert and Weiss studied a
special staircase surface, shown in Figure \ref{fig2} on the left \cite{HuWe}. They
showed that the straight-line flow is almost surely recurrent and completely classified the ergodic measures as well as
the periodic points.
In \cite{HoWe}, Hooper and Weiss classified the periodic square tiled staircases which are almost surely recurrent.  

In this article we study non-periodic staircases. 
We show that almost all staircases are recurrent. Two different notions of almost every staircase will be given, one for
square tiled staircases, the other more general.  We also show that the square tiled ones have dense set of directions
for which all regular orbits are periodic.   These results follow from approximating arbitrary staircases by periodic ones.
\section{Arbitrary staircases}
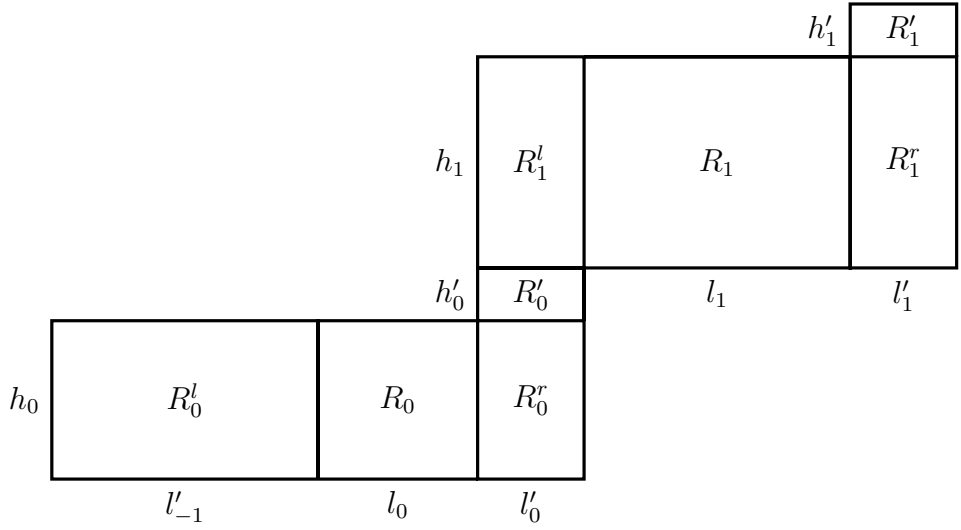
\begin{figure}[ht]
\centering
\begin{tikzpicture}[scale=0.7]
  \draw [very thick] (0,0)
    -- node [below] {$l_0$} ++ (-3,0)
    -- node [below] {$l'_{-1}$} ++ (-5,0)
    -- node [left] {$h_0$} ++ (0,3)
    --node [left] {} ++ (5,0)
    -- node [left] {} ++ (0,-3)
    -- node [left] {} ++ (0,3)
    -- node [below] {} ++ (3,0)
    -- node [below] {} ++ ( 0,-3 )
    -- node [below] {$l'_0$} ++ (2,0)
    -- node [below] {} ++ (0,3)
    -- node [below] {} ++ (-2,0)
    -- node [left] {$h'_0$} ++ (0,1)
    -- node [left] {} ++ (2,0)
    -- node [left] {} ++ (0,-1)
    -- node [left] {} ++ (0,1)
    -- node [left] {} ++ (-2,0)
    -- node [left] {$h_1$} ++ (0,4)
    -- node [below] {} ++ (2,0)
    -- node [below] {} ++ (0,-4)
    -- node [below] {$l_1$} ++ (5,0)
    -- node [below] {} ++ (0,4)
     -- node [left] {} ++ (-5,0)
      -- node [left] {} ++ (5,0)
    -- node [below] {} ++ (2,0)
    -- node [left] {} ++  (0,-4)
     -- node [below] {$l'_1$} ++ (-2,0)
     -- node [left] {} ++ (0,4)
    -- node [left] {$h'_1$} ++ (0,1)
    -- node [left] {} ++ (2,0)
    -- node [left] {} ++ (0,-1);
      \draw (-1.5,1.5) node {$R_0$};
      \draw (1,1.5) node {$R^r_0$};
      \draw (-5.5,1.5) node {$R^l_0$};
      \draw  (1,3.5) node {$R'_0$};
       \draw (1,6) node {$R^l_1$};
       \draw (4.5,6) node {$R_1$};
       \draw (8, 6) node {$R^r_1$};
       \draw (8, 8.5) node {$R'_1$};
\end{tikzpicture}
\caption{An arbitrary staircase}\label{fig1}
\end{figure} Consider $\Sigma' := (\R_+ \times \R^*_+ \times \R^*_+ \times \R_+)^\Z$ with the product topology (here $\R_+ :=\{x: x \ge0\}$ and $\R^*_+  := R_+ \setminus \{0\}$). Note that $\Sigma'$ is not compact. 
Fix $v := (l,h,l',h') \in \Sigma'$.
The staircase $T_v$ will be formed as follows (see Figure \ref{fig1}). All rectangles are oriented to have 
horizontal and vertical sides.  
There are four types of rectangles: the rectangle $R_n$ has length $l_n$ and height $h_n$,  $R'_n$ has dimensions $l'_n,h'_n$,
$R^l_n$ has
dimensions  $l_n',h_n$,  and $R^r_n$ has dimensions $l'_{n-1},h_n$.
The staircase consists of, for each $n \in \Z$, one copy of each of these rectangles.
We place the rectangle $R^r_{n}$ to the right of rectangle $R_{n}$  and we place $R^l_n$ to the left of $R_n$. Above $R^r_n$ we
place $R'_n$ while we place  $R'_{n-1}$ below $R^l_n$.
Continue this procedure inductively.
Touching edges and
edges which are a horizontal or vertical translation of one another
are identified.

We consider the straight-line flow $\psi_t$ on $T_v$  in a fixed direction $\theta \in \S1$. The phase volume is
the natural invariant measure of this flow.
A staircase is called {\em recurrent in the direction $\theta$} if for any set of positive measure in the phase space
a.e.\ orbit returns to this set and {\em recurrent} if it is recurrent in almost every direction.

Consider the shift transformation $\sigma: \Sigma' \circlearrowleft$ defined by $\sigma(v)_n = v_{n+1}$. 
Consider any $\sigma$-invariant ergodic probability measure $\nu$ on $\Sigma'$ whose support contains at least one
periodic staircase.  For example, let $f: \R_+ \circlearrowleft$ be a continuous Lebesgue integrable function with
integral 1 and let $\lambda$ be the probability measure on $\R_+$ whose Radon Nykodym derivative with respect to the length
measure is $f$.  Then 
$\lambda^{(4)} := \lambda  \times \lambda \times \lambda \times \lambda$ is a probability measure on $\R_+ \times \R^*_+ \times \R^*_+ \times \R_+$,
and thus the product measure $\nu_f$ determined by $\lambda^{(4)}$ is
a natural shift invariant probability measure on staircases. Since this measure is Bernoulli, periodic points are dense in its support.

\begin{theorem}\label{main1}
For $\nu$ a.e.~$v \in \Sigma'$ the staircase $T_v$ is recurrent.  
\end{theorem}
Remarque: we will show the the recurrent directions of the periodic staircase in the support of $\nu$ are recurrent directions for $T_v$.

\section{Square tiles staircases}
We consider a special case of the above construction with a different coding.
Let $\Si := \{0,1\}^{\Z}$ and endow $\Si$ with the product topology.
For each  $s \in \Si$ we construct the 
flat surface $T_s$ as follows  (see Figure 2). The surface $T_s$ consists of an infinite number of unit squares labelled by $\Z$.
The sides of the squares are parallel to the $x$ and $y$ axes.
The $n+1$st square is placed above the $n$th square if
$s_n = 0$ and otherwise it is to the right of the $n$th square.
The surface $T_s$ is obtained by identifying the common edges and identifying  pairs of edges which
are a horizontal or vertical translation
of one another.
For each $s \in \Sigma_2$ other than those ending with an infinite sequence of $0$'s or $1$'s (in the positive or negative direction) 
we can find a unique $v(s) \in \Sigma'$ which describes the same staircase with
with  $h_n = l'_n = 1$ and $h'_n \in \N, l_n \in \N$ for all $n$.  We leave the exact computation to the reader.
\begin{figure}[ht]
\begin{minipage}[ht]{0.49\linewidth}
\centering
\begin{tikzpicture}[scale=0.7]
  \draw [very thick] (0,0)
    -- node [below] {} ++ (-2,0)
    -- node [left] {} ++ (0,2)
    -- node [below] {} ++ (2,0)
    -- node [below] {} ++ ( 0,-2 )
    -- node [below] {} ++ (2,0)
    -- node [below] {} ++ (0,2)
    -- node [below] {} ++ (-2,0)
    -- node [left] {} ++ (0,2)
    -- node [below] {} ++ (2,0)
    -- node [below] {} ++ (0,-2)
    -- node [below] {} ++ (2,0)
    -- node [below] {} ++ (0,2)
    -- node [below] {} ++ (-2,0)
    -- node [left] {} ++ (0,2)
    -- node [left] {} ++ (2,0)
    -- node [left] {} ++ (0,-2);
      \draw (-1,1) node {$1$};
       \draw (1,1) node {$0$};
       \draw (1,3) node {$1$};
        \draw (3,3) node {$0$};
       \end{tikzpicture}
\end{minipage}
\begin{minipage}[ht]{0.49\linewidth}
\centering
\begin{tikzpicture}[scale=0.7]
  \draw [very thick] (0,0)
    -- node [below] {} ++ (-2,0)
    -- node [left] {} ++ (0,2)
    -- node [below] {} ++ (2,0)
    -- node [below] {} ++ ( 0,-2 )
    -- node [below] {} ++ (2,0)
    -- node [below] {} ++ (0,2)
    -- node [below] {} ++ (-2,0)
    -- node [left] {} ++ (2,0)
    -- node [below] {} ++ (2,0)
    -- node [below] {} ++ (0,-2)
    -- node [below] {} ++ (-2,0)
    -- node [below] {} ++ (0,2)
    -- node [left] {} ++ (0,2)
    -- node [below] {} ++ (2,0)
    -- node [below] {} ++ (0,-2)
        -- node [below] {} ++ (2,0)
-- node [below] {} ++ (0,2)
-- node [below] {} ++ (-2,0)
-- node [left] {} ++ (0,2)
-- node [below] {} ++ (2,0)
-- node [below] {} ++ (0,-2);
 \draw (-1,1) node {$1$};
       \draw (1,1) node {$1$};
              \draw (3,1) node {$0$};
                     \draw (3,3) node {$1$};
                       \draw (5,3) node {$0$};

\end{tikzpicture}
\end{minipage}

\caption{\mbox{The staircases $\cdot\cdot1010\cdot \cdot$ and $\cdot\cdot11010\cdot\cdot$}}\label{fig2}\end{figure}

We consider the shift transformation $\sigma$ on $\Sigma_2$. 
Consider any $\sigma$-invariant ergodic probability measure $\mu$ on $\Sigma$ whose support contains a
periodic staircase other than $0^\infty$ or $1^\infty$. 
Any measure with full support, in particular the Bernoulli measures, satisfy this
assumption.
\begin{theorem}\label{main2}
For $\mu$ a.e.~$s \in \Sigma_2$ the staircase $T_s$ is recurrent.
\end{theorem}

A direction $\theta \in \mathbb{S}^1$ is called {\em purely periodic} if all regular orbits in this direction are periodic. 
We additionally assume that the measure $\mu$ includes the periodic point $(10)^\infty$ (satisfied by Bernoulli measures).
\begin{theorem}\label{main3}
For  $\mu$ a.e.~$s \in \Sigma_2$, the staircase $T_s$ has a dense set of purely periodic directions. 
\end{theorem}
Remark: we will show that the set of purely periodic directions of $T_s$ includes the set of purely periodic directions of $T_{(10)^\infty}$
which was shown to be dense by \cite{HuWe}.

\section{Proofs}
The $\omega$ limit set of $v \in \Sigma'$
is $\cap_{n \ge 0} \overline{\{\sigma^kv: k > n\}}$ and the $\alpha$ limit set is defined similarly with $\sigma^k$ replaced by $\sigma^{-k}$.  Theorems \ref{main1} and \ref{main2} will follow almost immediately from the following result.

\begin{theorem}\label{recur2}
If there is a periodic point $v^+$ in the intersection of the  $\alpha$ and $\omega$ limit sets of $v \in \Sigma'$,  then the staircase $T_v$ is recurrent for any direction $\theta$ for which $T_{v^+}$ is recurrent.
\end{theorem}

\begin{proofof}{Theorem \ref{main1}}
By standard arguments
on skew products  due to Schmidt (\cite{Sc}, Theorem 11.4) (see also Proposition 10 of \cite{HoWe}) which we will sketch 
below, the staircase $T_{v^{+}}$ is recurrent in a.e.~direction $\theta$ , and thus Theorem \ref{main1} follows from Theorem \ref{recur2}.

Suppose $(X,\mu)$ is a finite measure space and $R: X \to X$ is a measurable transformation preserving $\mu$ which is ergodic.   For a measurable $f: X \to \Z$, $f \in L^1(X,\mu)$ define $X_f := X \times \Z$ and 
$$R_f : X_f \to X_f, \quad R_f(x,k) = (Rx,k+ f(x)).$$
Schmidt showed that $R_f$ is recurrent if and only if $\int f\, d\mu = 0$. We apply this result to a periodic staircase in the following way.  Consider a finite
crossection, say the side of one of our rectangles.  If we quotient the staircase $T_{v^{+}}$ by its period, then the first return map (for a fixed direction)  to this crosssection is an interval exchange map $R$, and Kerckhoff, Masur and Smillie have shown that
this interval exchange map is ergodic
for almost every direction  \cite{KMS}.  If we restrict to an ergodic direction, and consider the lifts of this cross section to 
$T_{v^{+}}$ we get a map of the form $R_f$.
From  our staircase construction we have $\int f \, d\mu = 0$ (this is holonomy 0 in the language of \cite{HoWe}), thus we
can apply Schmidt's theorem to conclude the recurrence of periodic staircases.
\end{proofof}

\begin{proofof}{Theorem \ref{main2}}
The $\omega$ (resp.\ $\alpha$) limit set of $s \in \Sigma_2$ includes a periodic point other that $0^\infty$ and $1^\infty$
if and only if the $\omega$ (resp.\ $\alpha$) limit set of the corresponding $v(s) \in \Sigma'$ contains a periodic point. Thus we can
conclude the proof by applying
Theorem \ref{recur2}.
\end{proofof}

\begin{proofof}{Theorem \ref{recur2}}
Let  $T = T_v$ and $T^+ = T_{v^+}$.
Fix a recurrent direction $\theta$ for $T_{v^+}$.
For each $n \in \Z$, let $L_n$  
denote the left boundary of the rectangle $R^r_n \cup R'_n \cup R^l_{n+1}$.  
Let  $\L_n = L_n \times \theta$ and let $\phi$ be the first return map of the straight-line 
flow on $T$ to the set $\L := \cup_{n \in \Z} \L_n$.
Denote corresponding objects in $T^+$ with a superscript  $+$, for example $\L^+$ and $\phi^+$. By 
assumption the map $\phi^+$ is recurrent.

The idea of the proof is quite simple. First we make quantitative estimates on the recurrence of $\phi^+$ on the periodic table $T^+$.
There are arbitrarily large finite parts of $T$ which are $\e$ close to a part of $T^+$, thus we can transfer
these estimates to $T$.  Care must be taken for orbits which come to close to a corner, and thus
can hit a different sequence of sides in $T$ and $T^+$. Conceptually the transfer should be clear, but the formal proof 
is a bit technical since the measure spaces are not identical, but only close.
Since the map  is invertible this will imply its recurrence, and finally we give an argument
to conclude the recurrence of the flow.

Fix $\e > 0$. Fix $N  > 0$ so that the $\phi^+$ orbit of at least $(1 - \e)\%$ of the points in $\L^+_0$ recurs to $\L^+_0$ before a time $N$.   Consider $D^+_N := \cup_{-N \le n \le N} L^+_n$ and $\mathcal{D}^+_N = D^+_N \times \theta$.  
There are $4$ corners of squares which are in the set $L^+_n$.  After identification, all corners of the table $T^+$
belong to the set $L^+$.  In total, there are $4(2N+1)$ corners in $D^+_N$ (without identifications).

Take a small $\e' = \e(h_0 + h'_0 + h_1)/(N\cdot 4(2N + 1))$ neighborhood $B_{\e'}$ of these corners in $D^+_N$. 
This set has measure (length) at most  $ \e(h_0 + h'_0 + h_1)/N$. Let  $\mathcal{B}_{\e'} := B_{\e'} \times \theta$.
Let
$C := \L_0 \cap \cup_{i=0}^{N-1} (\phi^+)^{-i} \mathcal{B}_{\e'} $, since $\phi^+$ is measure preserving, the measure of $C$ 
is at most $N \times \e(h_0 + h'_0 + h_1)/N  = \e(h_0 + h'_0 + h_1)$.  Since the total measure of $\L_0$ is
$h_0 + h'_0 + h_1$ we conclude that $(1-\e)\%$ of the points in $\L_0$ are not in $C$, 
i.e.~the first $N$ images of this $(1-\e)\%$ set of points
stay at least a distance $\e'$ away from 
the singular points.

Combining these two facts yields a set $G^+ \subset \L_0$ consisting of at least $(1 - 2\e)\%$ of the points in $\L^+_{0}$, such that the points of $G^+$ recur before time $N$ without having visited an $\e'$ neighborhood of a singular point.

Now consider the staircase $T$.  Let $H^+_n :=  h^+_n + {h'}^+_n +h^+_{n+1}$ be the length of $L^+_{n}$ (denote by $H_n$ the
corresponding length in $T$) and 
$H^+ = \inf_{n \in \Z} H_n $. Note that since $T^+$ is periodic $H^+$ is strictly positive.
Since $v^+$ is a $\omega$ limit point of $v$, for each $\e' > 0$ and $N > 0$ 
there is a positive $M$ so that 
$\max(|l_{i+M} - l_i^+|, |h_{i+M} - h_i^+| |,l'_{i+M} - {l'}_i^+|, |h'_{i+M} - {h'}_i^+|) <  \e'H^+/3N^2$ 
for all $|i| \le  N$.   Thus, for each $|i| \le N$ we have $|H^+_i - H_{M+i}| \le \e'H^+/N^2$.  For concreteness, fix the bottom point of each interval $L_{M+i}$ and $L^+_i$ as its origin.  Formally, this yields an
identification of $((H_{M+i}- \e'H^+/N)/H_{M+i})\% = (1 - \e'H^+/NH_{M+i})\%$ of $L_{M+i}$ with $L^+_i$ via the ``identity map'' 
for each fixed $|i| \le N$. In $\mathbb{R}^2$ this identification moves points by at most a distance of $N(\e' H^+/N^2) = \e' H^+/N$.

Consider a point $s^+ \in L^+_0$,  $s =id(s^+) \in L_M$ for which the identification exists.
Consider the $\phi$ and $\phi^+$ orbits of $(s,\theta)$.
Let $ (s_i,\theta_i) := \phi^i(s,\theta)$ and $(s^+_i,\theta^+_i) := (\phi^+)^i(s,\theta) $.
We compare these orbits for times $i \in \{0,1,\dots, N\}$.
As long as these orbits visit the ``same'' sequence of sides they remain parallel, i.e.~$\theta^+_i = \theta_i$
and the orbits diverge
at most linearly i.e.\ the distance $d(s^+_i,s_i) \le i(\e'H/N).$  If $0 \le i \le N$ then
$d(s^+_i,s_i) \le \e'H$.

In particular,  if the point $(s^+,\theta) \in G^+$ then for $0 \le i \le N$ the $s_i$ are a distance at least $\e'$
from the corners, so since $d(s_i,s_i) < \e'$ they are on the same side and we can apply this observation to
conclude that since $(s^+,\theta)$ is $\phi^+$-recurrent, the point $(s,\theta)$ is $\phi$-recurrent.  
This holds for at least $(1-3\e)\%$ of the
points in $\L_M$

Consider $B := \{(s,\theta) \in \L_0: s_i \not \in L_0 \text{ for all } i > 0 \}$, i.e.\ the set of points in $\L_0$ which do not
recur to $\L_0$. Since $\phi$ is invertible, we have $\phi^i B\cap \phi^jB = \emptyset$ for all $i > j$ (otherwise $\phi^{i-j}B \cap B \ne \emptyset$,
which contradicts the definition of $B$).  Fix $n > 0$. Thus a.e.\ point of $B$ can hit each level $\L_n$ $(n \ge 0)$ only a finite number of times.
Let $m_n(s,\theta)$ be the last time the orbit of $(s,\theta)$ visits $\L_n$, i.e.\ $s_{m_n} \in L_n$ and $s_j \not \in L_n$ for all $j > m_n$.
Then $\Phi(s,\theta) := \phi^{m_n(s,\theta)}(s,\theta)$ is a measure preserving map $\Phi : B \to \L_n$ defined almost surely. By definition
the image $\Phi(B)$
consists of nonrecurrent points in $\L_n$.  In particular if $n = M$ as chosen above then we conclude that the measure of $B$ is at most
$3 \e$.
Since $\e > 0$ was arbitrary, a.e.~point in $\L_0$ is $\phi$-recurrent.  We can repeat this argument with $\L_0$ 
replaced by $\L_j$ for any fixed $j$ to conclude that a.e.~point in $\L_j$ is $\phi$-recurrent for all $j$. 

Finally we need to show that the flow $\psi_t$ is recurrent in the direction $\theta$.
Consider any small open ball $B$ in the phase space of the flow $\psi_t$.  Flow each
non-singular point of  $B \times \{\theta\}$ 
until it hits the set $\L$. Since the ball is open, it has positive phase volume, and its image 
on the set $\L$ has positive phase area.   Almost every of these points is a $\phi$-recurrent point by the above. 

Fix a nonsingular $x \in B$ such that $(x_N,\theta) := \psi_t(x,\theta) \in \L_N$.  Note that
by transversality and the Fubini theorem almost every $x \in B$ corresponds to a $\phi$-recurrent $(x_N,\theta) \in \L_N$.  To conclude the proof we suppose that
$(x_N,\theta)$ is $\phi$-recurrent and we will show that this implies that  $(x,\theta)$ is $\psi_t$-recurrent.
Choose a  open neighborhood $U$ of $x$ small enough that for each $y \in U$ there is a $t(y)$ very close to $t$ such that
$\psi_{t(y)}(y,\theta)  \in \L_N$.  Let $U' = \{\psi_{t(y)}(y,\theta) : y \in U\}$.  This is a small neighborhood of $\psi_t(x,\theta)$, and by the above
results there is an (arbitrarily large) $n$ such that $\phi^n (x_N,\theta)  \in U'$. Thus $\phi^n (x_N,\theta) = \psi_s(x_N,\theta) = \psi_{s+t}(x,\theta) \in U'$ for some large $s$.   
Since  $\psi_{s+t}(x,\theta)$ is  in $U'$ it is the image $\psi_{t(y_0)}(y_0,\theta)$ for some $y_0 \in U$.  Thus  
$\psi_{s+t - t(y_0)}(x,\theta) = (y,\theta)$ with $y \in U$
and we conclude that $(x,\theta)$ is $\psi$ recurrent.
\end{proofof}

A purely periodic direction is called {\em strongly parabolic} if the phase space 
decomposes into
an infinite number of cylinders isometric to each other.  A periodic staircase is called good if there are a dense set
of  strongly parabolic directions.

\begin{theorem}\label{periodic}
If the $\alpha$ and $\omega$ limit sets of $s \in \Si$ include a periodic point $\hat{s} \in \Si$ such that 
the staircase $T_{\hat{s}}$ has a dense set of strongly periodic directions,
then the staircase $T_s$ has a dense set of purely periodic directions.
\end{theorem}

\begin{proofof}{Theorem \ref{main3}}
The theorem follows from Theorem \ref{periodic} since Hubert and Weiss have shown that the table $(10)^\infty$ satisfies the assumptions of Theorem \ref{periodic}.
\end{proofof}

\begin{proofof}{Theorem \ref{periodic}}
We consider another cross section to the billiard flow, let $\D_n$ be the left side of the $n$th square of the table $T_s$
and let $\D = \cup_{n \in \Z} \D_n$. Similarly let $\hat{\D}_n$ be the left side of the $n$th square of $\Ts$.  Clearly all non-singular vertical orbits are periodic in $T_s$.   Fix a strongly parabolic direction $\theta \in \S1$ for $\Ts$ which is not vertical, 
the flow in the direction $\theta$ in $T_s$ is transverse to $\D$.
Consider the first return map $\zeta: \D \times \{\theta\} \to \D \times \{\theta\}$ of the flow in this direction.  It is of the form
$$\zeta(r,n,\theta) = (r + \alpha \hspace{-10pt} \pmod 1, n + f(s,n,\theta),\theta)$$
for a certain $\alpha = \alpha(\theta)$.   In the table
$T_s$ one sees arbitrarily large initial pieces of the staircase $\Ts$: for each $N > 0$ there is are $M^- < 0 < M^+$ 
such that $s_{i+M^{\pm}} = \hat{s}_i$ for each $i \in \{-N,\dots,+N\}$. We immediately conclude that since the direction $\theta$ is strongly parabolic
for $\Ts$ the staircase $T_s$ must have periodic orbits in the direction $\theta$ and thus 
the number $\alpha$ must be rational.  

Since $\theta$ is strongly parabolic for $\Ts$ there is a $t_0$ so that all orbits in the direction $\theta$ have
flow period $t_0$.  Thus we can find an $N$ so that all nonsingular orbits starting in $\hat{\D}_0$ are periodic 
do not visit any $\hat{\D}_i$ with $|i| > N$. Thus we can conclude that for the table 
$T_s$ all nonsingular orbits starting the the two squares $T_{M^{\pm}}$ are periodic.

Now consider any point $x=(r,i,\theta)$ with nonsingular orbit with $M^- < i < M^+$.
First of all since $\D_{M^{\pm}}$ consists completely of
periodic orbits the orbit of $x$ can not reach $\D_{M^{\pm}}$ without being one of these periodic orbits.  If the orbit does not reach 
$\D_{M^{\pm}}$ then it 
stays inside the compact region $\cup_{i = M^- + 1}^{M^+ -1} \D_i$.  Since $\alpha$ is rational, and the number of $\D_i$
it visits is finite,  the  coordinates $(r_i,n_i)$ of the orbit 
can only take a finite number of values.  Thus it must visit some
point $(r,n)$ twice. Since the dynamics is invertible the orbit is periodic.
The result follows since $M^- \to - \infty$ and $M^+ \to \infty$ as $N \to \infty$
\end{proofof}

\section{Remarks}

These proofs hold in a more general setting,  they use only the one dimensionality of the construction of staircases. 
For the sake of clarity of exposition we did not try to define a larger class of flat surfaces.

The main idea of proof of this article 
was developed in a two dimensional setting for the Ehrenfest wind-tree model in \cite{Tr1}.  Since the approximation
technique is essentially one dimensional in the two dimensional setting it yields only topological results, i.e. a dense $G_{\delta}$ of recurrent
wind-tree models with dense purely periodic directions.  The recurrence result also holds for a dense $G_{\delta}$ of Lorentz gases \cite{Tr1}. 
 
Returning
to the one dimensional setting we see that our arguments also yield recurrence for almost every quench Lorentz tube
(see \cite{CrLeSe}).  The typical quenched Lorentz tubes do not needed to have finite horizon unlike those of \cite{CrLeSe}.
Using the hyperbolic structure, it can be shown that these Lorentz tubes are not only recurrent, but also ergodic \cite{LeTr}.

\end{document}